\documentclass{article}
\usepackage{amsthm}
\usepackage{amsfonts}
\usepackage{cite}
\usepackage{amsmath, amscd}
\newtheorem{theorem}{Theorem}[section]

\newtheorem{proposition}[theorem]{Proposition}
\newtheorem{lemma}[theorem]{Lemma}

\newtheorem{question}[theorem]{Question}

\begin{document}
\author{Jeremy Berquist}
\title{Demi-Normal Surface Singularities}
\maketitle

\noindent
\linebreak
\textbf{Abstract.  }   We prove semi-rationalification and semi-log-canonicalization for Gorenstein demi-normal surfaces.  That is, given a Gorenstein demi-normal surface $X$ with semi-rational (respectively, semi-log canonical) singularities in an open set $U$ with complement a finite set of points, there is a proper birational morphism $f:Y \rightarrow X$  such that $f$ is an isomorphism over $U$ and $Y$ has only semi-rational (respecitvely, semi-log canonical) singularities. .  We proceed by passing to the normalization and then gluing along the conductor in an appropriate rationalification or log-canonicalization of $\overline{X}$.  It is not hard to prove the analogous results in the normal surface case using Zariski's method of resolving normal surface singularities (that is, by repeatedly blowing up points and then normalizing, a process that halts after finitely many steps with a smooth surface), and by using the log canonicalization of Hacon and Xu.  We prove that a demi-normal variety $X$ has semi-rational singularities if its normalization $\overline{X}$ has rational singularities, and that the converse is true when $X$ is Gorenstein.  Likewise, we prove that $X$ is semi-log canonical if and only if the pair $(\overline{X}, C)$ is log canonical, where $C$ is the conductor, and $X$ has  $\mathbb{Q}$-Cartier canonical class.  We use the following terminology.  A demi-normal surface has semi-rational singularities if there is a semi-resolution $f: Y \rightarrow X$ such that $R^if_*\mathcal{O}_Y = 0$ for $i>0$; equivalently, since $X$ is assumed to be $S_2$ (hence Cohen-Macaulay), $f_*\omega_Y = \omega_X$.  A semi-smooth Grauert-Riemenschneider vanishing result is required here.  We say that $X$ is semi-log canonical if $K_Y = f^*K_X + \Sigma a_iE_i$ with $a_i \geq -1$ for all $i$, where $E_i$ are the exceptional divisors of $f$.  The Gorenstein hypothesis amounts to the fact that $\omega_X$ is an invertible sheaf, and hence the canonical class $K_X$ is Cartier.  Finally, we give unproved conditions under which hold the analogous results for varieties in higher dimensions.

\newpage
\tableofcontents
\addcontentsline{}{}{}
\newpage
\begin{section}{Introduction}
We study quasi-projective demi-normal surfaces.  These are taken as reduced varieties of dimension two over a field of characteristic zero, with Serre's $S_2$ property, and such that there is an open subvariety $U$ with complement of codimension two (i.e., a finite set of closed points) such that the only singularities in $U$ are analytically double normal crossings.  Equivalently, these are surface singularities which are $S_2$, $G_1$ (Gorenstein in codimension one), and SN (seminormal).  We assume also that such a surface $X$ is Gorenstein. That is, that its dualizing sheaf $\omega_X$ (which exists, since $X$ is Cohen-Macaulay by definition) is invertible.   It is known that both semi-rational and semi-log canonical surfaces are $\mathbb{Q}$-Gorenstein, meaning that some tensor power of $\omega_X$ is invertible, so our assumption amounts to the fact that the index is 1.  We might then replace $X$ with its index-1 cover, but such a procedure is not birational, so we seem to need this stronger condition at the outset.  

A demi-normal variety has semi-rational singularities if for a given semi-resolution (and hence for every semi-resolution) $f: Y \rightarrow X$, we have $R^if_*\mathcal{O}_Y = 0$ for $i>0$.  Note that $f_*\mathcal{O}_Y =\mathcal{O}_X$ for every semi-resolution.  An equivalent definition is that $X$ is Cohen-Macaulay and $R^if_*\omega_Y = 0$ for $i>0$ and $f_*\omega_Y = \omega_X$.  The first of these conditions, Grauert-Riemenschneider vanishing, holds automatically.  See \cite{Berq14}.  Thus a demi-normal surface is semi-rational if $f_*\omega_Y = \omega_X$.  We say that a normal surface is rational if the same statements are true for a resolution of singularities.  Again, it is well-known that GR-vanishing is automatic.  See \cite{GR70}.

A demi-normal variety is semi-log canonical if, given a semi-resolution as above, we have $K_Y = f^*K_X + \Sigma a_iE_i$, where $a_i \geq -1$ for all $i$.  This is again independent of the semi-resolution chosen.  As part of the definition, one requires that $X$ is $\mathbb{Q}$-Gorenstein, so that $f^*K_X = \frac{1}{m}f^*(mK_X)$ is obtained by pulling back sections of the Cartier divisor $mK_X$.  The definition of log canonical singularities is the same, except that we use a resolution of singularities of a normal variety, instead of a semi-resolution of a demi-normal variety.

Zariski's method for resolving the singularities of a normal surface is to repeatedly blow up points and then normalize.  He shows that this procedure results in a smooth surface after finitely many iterations.  See \cite{Zar39}.  Thus it is easy to show that a normal surface has a rationalification.  Let $U$ be the open subvariety of $X$ on which $X$ has rational singularities.  The complement of $U$ is necessarily a finite set of points.  By repeatedly blowing up points and then normalizing, we do not alter $U$, and the process stops when the non-rational singularities of $X$ are replaced by smooth points.  In other words, if the singular set of $X$ is the set of points $\{ p_1, \ldots, p_s, p_{s+1}, \ldots, p_n \}$, where $p_1, \ldots, p_s$ are rational, we apply Zariski's procedure to $X-\{p_1, \ldots, p_s\}$.  

Log canonicalization in all dimensions has been proved by Hacon and Xu.  We use their result for the surface case in order to construct a semi-log canonicalization and explain why the analogous procedure of passing to the normalization and gluing along the conductor does not work for dimensions greater than two.

There is a good relation between rational and semi-rational singularities for all dimensions.  Namely, if the normalization $\overline{X}$ has rational singularities, and both $X$ and $\overline{X}$ are Cohen-Macaulay, then $X$ has semi-rational singularities.  We prove this in the next section.  Thus for demi-normal surfaces, where we obtain the Cohen-Macaulay condition for free, this allows us to look at the canonical sheaves alone.  However, the converse is not necessarily true.  It is here where we use the Gorenstein condition on $X$.  If $X$ has Gorenstein semi-rational singularities, then it has semi-canonical sigularities.  Thus the pair $(\overline{X}, C)$ is canonical.  In particular, it is log terminal, and such singularities are known to be rational.  See the lecture notes by Koll\'ar, ``Singularities of Pairs," Theorem 11.1.

We also have a correspondence between log canonical and semi-log canonical singularities.  Namely, when $X$ is semi-log canonical, the pair $(\overline{X},C)$ is log canonical, and vice-versa, provided that we can insure $X$ is $\mathbb{Q}$-Gorenstein.  

The method of proof for both procedures on $X$ is to pass to the normalization.  We use rationalification and log-canonicalization of normal surfaces, as described above.  Then we glue along the conductor in order to obtain a proper birational morphism $f: Y \rightarrow X$ such that the normalization $\overline{Y}$ of $Y$ has either rational or log canonical singularities.  Finally, we use the correspondence between singularities on $Y$ and $\overline{Y}$ in order to show that $Y$ has the right singularities.  The tool we need is a gluing procedure.  Given a variety $Z$ with a closed subvariety $W$, and a finite morphism $W \rightarrow W_0$, there is a universal pushout 
$$\begin{CD}
W      @>>>  Z \\
@VVV             @VVV \\
W_0  @>>> Z_0 \\
\end{CD}$$
such that $Z \rightarrow Z_0$ is finite, agrees with $W \rightarrow W_0$ on $W$, and is an isomorphism away from $W$.  Here we take $Z$ to be a normal  surface with specified (rational or log canonical) singularities and $W$ to be a birational transform of the conductor $C$ in $\overline{X}$, where $X$ is the demi-normal variety we begin with.   In particular, $Z$ is the normalization of $Z_0$.  We may even replace $Z_0$ by its demi-normalization to assume that it is demi-normal.  See \cite{Berq14b}.  For facts about the gluing construction, see \cite{Art70}.

The problem in higher dimensions is that there is no finite morphism $W \rightarrow W_0$ to glue along.  For surfaces, where $W$ is a one-dimensional variety, there is no issue.  A proper morphism with finite fibers is finite.  But if the dimension is larger than two, then the available morphism $W \rightarrow W_0$ is proper, but typically with some fibers that are at least one-dimensional, so the morphism is not finite.  This seemingly minor inconvenience is at the heart of our troubles with extending the surface case to higher dimensions.  We will discuss this problem in more detail in the final section, having seen how the surface case is rather special.

The paper is organized as follows.  In the second section, we prove semi-rationalification for Gorenstein demi-normal surfaces.  In the third section, we prove semi-log canonicalization for these varieties.  Finally, we outline a possible proof for the higher-dimensional case, citing both open problems (rationalification is not known to be true even for 3-folds) and known results (Hacon and Xu have proved log canonicalization in arbitrary dimensions, see \cite{HXu11}), as well as a possible procedure for proving semi-rationalification and semi-log canonicalization in general, given the results for normal varieties.
\end{section}

\begin{section}{Semi-Rationalification of Surfaces}

Let $X$ be a Gorenstein demi-normal surface with semi-rational singularities in an open set $U$.  Since $X$ is semi-smooth in codimension one, the complement of $U$ is a finite set of points.  Let $p:  \overline{X} \rightarrow X$ be the normalization of $X$, and let $C \hookrightarrow \overline{X}$ and $D \hookrightarrow X$ be the conductors.  We first prove that $p^{-1}U$ has rational singularities.  

\begin{proposition}  Suppose $X$ is a Gorenstein demi-normal surface with semi-rational singularities.  Then $\overline{X}$ has rational singularities.  The converse holds even without the Gorenstein hypothesis.

\proof  Let $f: Y \rightarrow X$ be a semi-resolution.  We have an induced commutative diagram
$$\begin{CD}
\overline{Y}  @>\overline{f}>>  \overline{X} \\
@VqVV                                           @VpVV \\
Y                   @>f>>                     X \\
\end{CD}$$.

Here $q$ is the normalization of a semi-smooth variety $Y$, and hence $\overline{Y}$ is smooth.  Since $\overline{f}$ is proper and birational (all other morphisms are birational, the composition in either direction is proper, and $p$ is separated), it is a resolution of singularities.

Now suppose that $X$ is Gorenstein with semi-rational singularities.  In particular, we have $f_*\omega_Y = \omega_X$.  Since $\omega_X$ is invertible, we have an injection $f^*\omega_X \rightarrow \omega_Y$.  Hence $X$ is in particular semi-canonical.  Then the pair $(\overline{X},C)$ is canonical (this will be proved in more detail in the next section).  In particular, the pair is log terminal, and it is well-known that this implies $\overline{X}$ has rational singularities.

For the converse, suppose $\overline{X}$ has rational singularities.  We prove that $X$ has semi-rational singularities.   Since $X$ is Cohen-Macaulay, it is enough to show that $f_*\omega_Y = \omega_X$.  We know that $\overline{f}_*\omega_{\overline{Y}} = \omega_{\overline{X}}$.  Let us follow $\omega_{\overline{Y}}$ in both directions in the commutative square above.

On the one hand, we have $$p_*\overline{f}_*\omega_{\overline{Y}} = p_*\omega_{\overline{X}} = \mathcal{H}om(p_*\mathcal{O}_{\overline{X}}, \omega_X),$$  using duality for a finite morphism.   In the opposite direction, we have $$f_*q_*\omega_{\overline{Y}} = f_*\mathcal{H}om(q_*\mathcal{O}_{\overline{Y}}, \omega_Y).$$  This last sheaf maps into $\mathcal{H}om(f_*q_*\mathcal{O}_{\overline{Y}}, f_*\omega_Y) = \mathcal{H}om(p_*\mathcal{O}_{\overline{X}}, f_*\omega_Y)$.  In particular, we obtain a morphism $$\mathcal{H}om(p_*\mathcal{O}_{\overline{X}}, \omega_X) \rightarrow \mathcal{H}om(p_*\mathcal{O}_{\overline{X}}, f_*\omega_Y).$$  Composing with the morphism induced by the natural inclusion $f_*\omega_Y \hookrightarrow \omega_X$, we obtain $$\mathcal{H}om(p_*\mathcal{O}_{\overline{X}}, \omega_X) \rightarrow \mathcal{H}om(p_*\mathcal{O}_{\overline{X}}, f_*\omega_Y) \rightarrow \mathcal{H}om(p_*\mathcal{O}_{\overline{X}}, \omega_X).$$  Now the sheaf on either end is reflexive, since $\omega_X$ is reflexive, and the composition is obviously an isomorphism in codimension one, since $f$ is a semi-resolution.  We conclude that the composition is the identity, and in particular that the second morphism is surjective.

Since either $\mathcal{H}om$ maps surjectively onto the second factor via evaluation at 1, we have that the natural injection $f_*\omega_Y \hookrightarrow \omega_X$ is also a surjection.  We conclude that $X$ has semi-rational singularities.  \qed

\end{proposition}

Note that the above proof makes sense in any dimension.  The surface case is special for another reason, as we will see shortly.

We have proved that $p^{-1}U$ has rational singularities.  As we discussed in the introduction, Zariski's method of resolving surface singularities allows us to conclude that $\overline{X}$ has a rationalification.  In other words, there is a proper, birational morphism $f:  Z \rightarrow \overline{X}$ such that $Z$ has only rational singularities and $f$ is an isomorphism over $p^{-1}U$.  Let $C'$ be the birational transform of the conductor $C$ in $\overline{X}$. Then the composition $C' \rightarrow C \rightarrow D$ is a finite morphism.  In fact, these are all curves, and the proper morphism $C' \rightarrow C$ has finite fibers.  Thus we may form the pushout 
$$\begin{CD}
C'  @>>> Z \\
@VVV       @VVV \\
D  @>>> X_0 \\
\end{CD}$$

By the universal property of the pushout, we obtain a morphism $X_0 \rightarrow X$.  Furthermore, since $C' \rightarrow D$ agrees with $C \rightarrow D$ on the open set $p^{-1}U$, we see that $X_0 \rightarrow X$ is an isomorphism over $U$.  We may further replace $X_0$ by its demi-normalization $X_1$ in $Z$.  See \cite{Berq14}.  We obtain a morphism $X_1 \rightarrow X$ that is an isomorphism over $U$, and such that the normalization $Z$ of $X_1$ has rational singularities.  We conclude by the proposition that $X_1$ has semi-rational singularities.

We have proved:

\begin{proposition}  Suppose $X$ is a Gorenstein demi-normal surface with semi-rational singularities in an open set $U$.  Then $X$ has a semi-rationalification.

\end{proposition}

Note that in general, $C' \rightarrow C$, which can be thought of as a sequence of blowups outside of $C \cap p^{-1}U$, will not in general be finite.  This is why the surface case is special.  If we blow up points on a curve, the morphism is finite.  So far, we have not found a way to resolve this issue in higher dimensions.
\end{section}

\begin{section}{Semi-Log Canonicalization of Surfaces}
We prove, using a theorem of Hacon and Xu, that a demi-normal surface $X$ with semi-log canonical singularities in an open set $U$ has a semi-log canonicalization.  That is, there exists a proper, birational morphism $X' \rightarrow X$ such that $X'$ has only semi-log canonical singularities and which is an isomorphism over $U$.

We need a preliminary lemma.  We use the notation from the previous section.

\begin{lemma}  $X$ is semi-log canonical if and only if it is $\mathbb{Q}$-Gorenstein and the pair $(\overline{X}, C)$ is log canonical.

\proof  Letting $f:  Y \rightarrow X$ be a semi-resolution and $\overline{f}:  \overline{Y}  \rightarrow \overline{X}$ be the induced resolution of singularities, we pull back from a relation $$K_Y \sim f^*K_X + \Sigma a_i E_i.$$  Note that for the normalization $p$, we have $p^*\omega_X = \omega_{\overline{X}}(C)$, and likewise for the normalization $q$.  Thus we have the equivalent relation $$K_{\overline{Y}} + C' = \overline{f}^*(K_{\overline{X}} + C) + \Sigma a_i q^{-1}E_i.$$  Here $C'$ is the birational transform of $C$, and it is easy to see that the $q^{-1}E_i$ are the exceptional divisors of $\overline{f}$.  Thus we have the lemma.  Note that for a semi-resolution, no component of the conductor is exceptional.  \qed

\end{lemma}

We state the theorem due to Hacon and Xu.  See \cite{HXu11}, Corollary 1.2.

\begin{proposition}  Let $(W,C)$ be a normal variety such that for some open set $U$, the pair $(U, C|_U)$ is log canonical.  Then there exists a proper, birational morphism $f:  W' \rightarrow W$ such that $(W',C')$ is log canonical, and $(W'|_{f^{-1}U}, C'|_{f^{-1}U})$ is isomorphic to $(U, C|_U)$.
\end{proposition}

Now let $X$ be a Gorenstein, demi-normal surface with semi-log canonical singularities in an open set $U$.  By (3.1), the pair $(\overline{X}, C)$ is log canonical in $p^{-1}U$.  By (3.2), there is a log canonicalization of this pair.  Let $C'$ be the birational transform of $C$.  As in the previous section, we know that $C' \rightarrow C \rightarrow D$ is finite, and so we can glue along this morphism.  After possibly taking the demi-normalization, we can conclude the following:  there is a proper, birational morphism $X_1 \rightarrow X$ that is an isomorphism over $U$, and such that the normalization pair $(Z, C')$ is log canonical.  We would like to conclude that $X_1$ is semi-log canonical.  It remains to show that $X_1$ is $\mathbb{Q}$-Gorenstein.

For this, we use results of Koll\'ar.  In particular, we have the following.  See the lecture notes of Koll\'ar, ``Semi Log Canonical Pairs," 2010.

There is a Galois involution $\tau:  \overline{C} \rightarrow \overline{C}$ on the normalization $\overline{C}$ of $C$.  Since $K_X$ is Cartier, the different $\textnormal{Diff}_{\overline{C}}0$ is $\tau$-invariant. This is Proposition 14 in the above article.   Note in our construction that the map $C' \rightarrow C$ is proper and birational (hence finite), and so the normalization of $C'$ is also the normalization of $C$.  We have constructed a pair $(Z,C')$ that is log canonical and such that $Z$ is the normalization of the demi-normal variety $X_1$.  By Theorem 17 of the same article, there is a codimension 3 subset outside of which $X_1$ is semi-log canonical.  Since we are dealing with surfaces, this implies directly that $X_1$ is semi-log canonical.  In particular, it is $\mathbb{Q}$-Gorenstein.  In other words, the $\tau$-invariance of the different follows from the fact that $X$ is Gorenstein, and nothing is changed with respect to $\tau$ when replacing $C$ with its birational transform $C'$.  

We have proved:

\begin{proposition}  A Gorenstein demi-normal surface with semi-log canonical singularities in an open set $U$ admits a semi-log canonicalization.
\end{proposition}

We note that the proof breaks down in dimensions greater than two.  For one thing, we do not always have a finite morphism to glue along, as we remarked in the last section.  Even if we did, we would only be able to conclude that $X_1$ has semi-log canonical singularities outside a codimension 3 subset.  

\end{section}
\begin{section}{Concluding Remarks}
The gluing problem is at the heart of our problems with extending the surface results to higher dimensions.  We state it as an open question.

\begin{question}  Given a finite morphism $C \rightarrow D$ such that $C$ is a dense open subset of $C'$, under what conditions is there a finite morphism $C' \rightarrow D'$ that agrees with $C \rightarrow D$ when restricted to $C$?

\end{question}
When we perform a rationalification or a log canonicalization, we replace the conductor $C$ in $\overline{X}$ with its birational transform $C'$.  We would like to glue along this new object in order to obtain a semi-rationalification or semi-log canonicalization of $X$.  As we saw in the previous section, we might be in good shape if $C'$ is normal, since in that case there exists an involution on $C'$.  In the semi-log canonicalization situation, we are still stuck, because we can only conclude that $X$ is semi-log canonical in codimension two.  However, for semi-rationalification, where we do not require a $\mathbb{Q}$-Gorenstein hypothesis, we could proceed as follows.  Supposing $C'$ is normal, glue along the quotient $C' \rightarrow C'/\tau$.  Upstairs we have only rational singularities, and then by (2.1) we would have semi-rational singularities downstairs.  Of course, we need to assume that both $X$ and $\overline{X}$ are Cohen-Macaulay in order to use the analog of (2.1), a condition that is automatic in the surface case.  Criteria for Cohen-Macaulayness have been investigated in \cite{Berq15}.

\end{section}
\nocite{Zar39}
\nocite{HXu11}
\nocite{Berq15}
\nocite{Art70}
\nocite{Berq14}
\nocite{GR70}
\nocite{Berq14b}

\bibliographystyle{plain}
\bibliography{paper}

\begin{thebibliography}{1}

\bibitem{Art70}
M.~Artin.
\newblock Algebraization of {F}ormal {M}oduli: I{I}. {E}xistence of
  {M}odifications.
\newblock {\em Annals of Mathematics}, 91(1):88--135, 1970.

\bibitem{Berq14}
J.~Berquist.
\newblock On {S}emi-{R}ational {S}ingularities, arxiv:1408.5835.
\newblock 2014.

\bibitem{Berq14b}
J.~Berquist.
\newblock Projective {C}losures of {D}emi-{N}ormal {V}arieties,
  arxiv:1411.2264.
\newblock 2014.

\bibitem{Berq15}
J.~Berquist.
\newblock Rational and {S}emi-{R}ational {S}ingularities, arxiv : 1503.06320v1.
\newblock 2015.

\bibitem{GR70}
H.~Grauert and O.~Riemenschneider.
\newblock Verschwindungssatze fur analytische {K}ohomologiegruppen auf
  komplexen {R}aumen.
\newblock {\em Inventiones Mathematicae}, 11:263--292, 1970.

\bibitem{HXu11}
C.~Hacon and C.Xu.
\newblock Existence of {L}og {C}anonical {C}losures, arxiv:1105.1169v3.
\newblock 2012.

\bibitem{Zar39}
O.Zariski.
\newblock The {R}eduction of the {S}ingularities of an {A}lgebraic {S}urface.
\newblock {\em Ann. of Math.}, 2:639--689, 1939.

\end{thebibliography}

\end{document}